# Some integral identities involving products of general solutions of Bessel's equation of integral order

S.K.H. Auluck[#]
Advanced Technology Systems Section
Bhabha Atomic Research Center, Mumbai-400085, India

Abstract

Spectral decomposition of dynamical equations using curl-eigenfunctions has been extensively used in fluid and plasma dynamics problems using their orthogonality and completeness properties for both linear and non-linear cases. Coefficients of such expansions are integrals over products of Bessel functions in problems involving cylindrical geometry. In this paper, certain identities involving products of two and three general solutions of Bessel's equation have been derived. Some of these identities have been useful in the study of Turner relaxation of annular magnetized plasma [S.K.H. Auluck, Phys. Plasmas, 16, 122504, 2009], where quadratic integral quantities such as helicity and total energy were expressed as algebraic functions of the arbitrary constants of the general solution of Bessel's equation, allowing their determination by a minimization procedure. Identities involving products of three solutions enable expanding a product of two solutions in a Fourier-Bessel series of single Bessel functions facilitating transformation of partial differential equations representing non-linear dynamics problems into time evolution equations by eliminating spatial dependences. These identities are required in an ongoing investigation of fluctuation-driven coherent effects in a nonlinear dynamical system.

# e-mail: auluck@barc.gov.in

# I. Introduction

Eigenfunctions of the curl operator, also known as Chandrasekhar-Kendall (C-K) functions [1] and their generalization to other geometries, can be used as basis functions to expand vector fields in spectral decomposition methods because of their orthogonality and completeness as well as because of simplifications introduced in handling electromagnetic fields and fluid vorticity [2]. This is usually confined to linear equations because of the limitations involved in determining the spectral coefficients of such expansion. In cylindrical geometry, C-K functions have Bessel functions describing radial modes, sine and cosine functions describing axial modes in a periodic boundary condition and the azimuthal modes. If the dynamical variable are expressed as Fourier-Bessel series in C-K functions (each C-K function representing an eigenmode), then their products in non-linear equations would be described as a double series in products of C-K functions. If a product of C-K functions itself can be expressed as Fourier Bessel series in C-K functions, then the spectral decomposition method can in principle be extended to non-linear equations as well, facilitating transformation of partial differential equations representing non-linear dynamics problems into time evolution equations by eliminating spatial dependences. This involves numerically evaluating and tabulating integrals involving products of three C-K functions. However, when large mode numbers, representing finer spatial scales, are involved, the integrands vary rapidly causing difficulties in numerical evaluation. Analytical expressions for such integrals would be very useful in dealing with large mode number situations enabling a deeper understanding of scaling relationships in the parameter space. If analytical methods are not available, a large computational effort required to generate a database of such integrals using arbitrary precision methods once for all would be quite useful.

In this paper, an attempt is made to derive some integral identities involving products of Bessel functions, encountered in simple problems with azimuthally symmetric geometry. Some of these identities have been useful in the study of Turner relaxation of annular magnetized plasma [3], where quadratic integral quantities such as helicity and total energy were expressed as algebraic functions of the arbitrary constants of the general solution of Bessel's equation, allowing their determination by a minimization procedure. A part of the following is included in the appendix of [3], which is re-derived below for the sake of completeness and continuity.

# II. Identities involving products of two functions:

The functions

$$Z_n(x) \equiv aJ_n(x) + bY_n(x) \qquad 1$$

where a and b are arbitrary constants, obey the following relations common to $J_n$ and $Y_n$, the Bessel functions of first and second kinds with integral order n.

$$\frac{2n}{x}Z_n(x) = Z_{n-1}(x) + Z_{n+1}(x) \qquad 2$$

$$2\frac{dZ_n(x)}{dx} = Z_{n-1}(x) - Z_{n+1}(x) \qquad 3$$

$$Z_{-n}(x) = (-1)^n Z_n(x) \qquad 4$$

If Up($\alpha$x) and Vq ($\beta$x) are two solutions of Bessel's equation of order p and q, then the following identity is known [4] :

$$\begin{aligned}&\int\left\{\left[\left(\alpha^2-\beta^2\right)x - \frac{p^2-q^2}{x}\right]U_p(\alpha x)V_q(\beta x)\right\}dx \\ &= \beta x U_p(\alpha x)V_{q-1}(\beta x) - \alpha x U_{p-1}(\alpha x)V_q(\beta x) + (p-q)U_p(\alpha x)V_q(\beta x)\end{aligned} \qquad 5$$

Putting $U \equiv Z_n(\alpha x);\ V \equiv \mathcal{Z}_n(\beta x)$, (where $Z$ and $\mathcal{Z}$ could in general have different values of a and b, a situation which might be encountered in some problems) and using 2 and 3

$$\begin{aligned}\int x dx\left[Z_n(\alpha x)\mathcal{Z}_n(\beta x)\right] &= \frac{1}{\left(\alpha^2-\beta^2\right)}\left\{\beta x Z_n(\alpha x)\mathcal{Z}_{n-1}(\beta x) - \alpha x Z_{n-1}(\alpha x)\mathcal{Z}_n(\beta x)\right\} \\ &= \frac{x}{\left(\alpha^2-\beta^2\right)}\left\{Z_n(\alpha x)\frac{d\mathcal{Z}_n(\beta x)}{dx} - \mathcal{Z}_n(\beta x)\frac{dZ_n(\alpha x)}{dx}\right\}\end{aligned} \qquad 6$$

In the case of a problem with annular cylindrical geometry [2], where the normalized inner radius corresponds to x=1, the normalized outer radius corresponds to x=A, (=outer radius /inner radius) and physical boundary conditions require that some of the fields (expressed as linear combination of Bessel functions for a given order n) vanish at the boundaries: $Z_n(\alpha A) = Z_n(\alpha) = \mathcal{Z}_n(\beta A) = \mathcal{Z}_n(\beta) = 0$ the definite integral over the domain becomes

$$\int_1^A x dx\left[Z_n(\alpha x)\mathcal{Z}_n(\beta x)\right] = \frac{1}{(2\alpha)}\left\{\left(A^2\right)\frac{d\mathcal{Z}_n(\beta A)}{dx}\frac{dZ_n(\alpha A)}{dx} - \frac{d\mathcal{Z}_n(\alpha)}{dx}\frac{dZ_n(\alpha)}{dx}\right\}\delta_{\alpha,\beta} \qquad 7$$

This orthogonality relation, which is a direct consequence of the boundary conditions, allows expressing quadratic integral forms, defined over an annular domain, involving physical fields, such as total energy, power dissipation, helicity, stress tensor in terms of the arbitrary coefficients a,b which define the linear combination in equation 1 and which can be determined at the end of the problem through an extremum principle.

For subsequent discussion, $Z$ and $Z$ are taken to be identical. Substitute $\beta=\alpha+\varepsilon$ and expand in first power of $\varepsilon$ and take the limit as $\varepsilon\rightarrow 0$.

$$\int x dx\left[Z_n(\alpha x)\right]^2=\frac{x^2}{2}\left[\left[Z_n(\alpha x)\right]^2-Z_{n+1}(\alpha x)Z_{n-1}(\alpha x)\right] \tag{8}$$

For n=0 and n=1

$$\int x dx\left[Z_0(\alpha x)\right]^2=\frac{x^2}{2}\left[\left[Z_0(\alpha x)\right]^2+\left[Z_1(\alpha x)\right]^2\right] \tag{9}$$

$$\int x dx Z_1(\alpha x)^2=\frac{x^2}{2}\left\{\left[Z_1(\alpha x)\right]^2+\left[Z_0(\alpha x)\right]^2-\frac{2}{\alpha x}Z_0(\alpha x)Z_1(\alpha x)\right\} \tag{10}$$

Multiplying 2 and 3 for $x\rightarrow\alpha x$, recurrence relation is obtained for higher order integrals similar to 9, 10:

$$\int x dx\left[Z_{n+1}(\alpha x)\right]^2=\int x dx\left[Z_{n-1}(\alpha x)\right]^2-\left[Z_n(\alpha x)\right]^2\frac{2n}{\alpha^2} \tag{11}$$

Using relations $\frac{d}{dx}Z_0(\alpha x)=-\alpha Z_1(\alpha x);\ \frac{d}{dx}xZ_1(\alpha x)=\alpha x Z_0(\alpha x)$ obtained from 2 and 3, the following can be derived

$$\begin{aligned}&\int x^p Z_1^2(\alpha x)dx\\&=-\frac{1}{\alpha}x^p Z_1(\alpha x)Z_0(\alpha x)+\frac{(p-1)}{2\alpha^2}\int x^{p-3}d\left[xZ_1(\alpha x)\right]^2+\int x^p\left[Z_0(\alpha x)\right]^2dx\end{aligned} \tag{12}$$

For p=3

$$\int x^3\left\{\left[Z_0(\alpha x)\right]^2-\left[Z_1(\alpha x)\right]^2\right\}dx=\alpha^{-1}x^3Z_1(\alpha x)Z_0(\alpha x)-\alpha^{-2}x^2\left[Z_1(\alpha x)\right]^2 \tag{13}$$

Using 8 with n=0 and 4,

$$\left[Z_0(\alpha x)\right]^2x=\frac{d}{dx}\left\{\frac{x^2}{2}\left[\left[Z_0(\alpha x)\right]^2+\left[Z_1(\alpha x)\right]^2\right]\right\} \tag{14}$$

Therefore

$$x^3\left[Z_0(\alpha x)\right]^2 = x^2\frac{d}{dx}\frac{x^2}{2}\left[\left[Z_0(\alpha x)\right]^2+\left[Z_1(\alpha x)\right]^2\right]$$
$$=\frac{d}{dx}x^2\frac{x^2}{2}\left[\left[Z_0(\alpha x)\right]^2+\left[Z_1(\alpha x)\right]^2\right]-x^3\left[\left[Z_0(\alpha x)\right]^2+\left[Z_1(\alpha x)\right]^2\right] \quad 15$$

Thus

$$\int x^3\left[2\left[Z_0(\alpha x)\right]^2+\left[Z_1(\alpha x)\right]^2\right]dx=\frac{x^4}{2}\left[\left[Z_0(\alpha x)\right]^2+\left[Z_1(\alpha x)\right]^2\right] \quad 16$$

Adding 16 and 13

$$\int x^3\left[Z_0(\alpha x)\right]^2 dx=\frac{x^4}{6}\left[\left[Z_0(\alpha x)\right]^2+\left[Z_1(\alpha x)\right]^2\right]+\frac{1}{3\alpha}x^3Z_1(\alpha x)Z_0(\alpha x)-\frac{1}{3\alpha^2}x^2\left[Z_1(\alpha x)\right]^2 \quad 17$$

Subtracting 13 from 17

$$\int x^3\left[Z_1(\alpha x)\right]^2 dx=\frac{x^4}{6}\left[\left[Z_0(\alpha x)\right]^2+\left[Z_1(\alpha x)\right]^2\right]-\frac{2}{3\alpha}x^3Z_1(\alpha x)Z_0(\alpha x)+\frac{2}{3\alpha^2}x^2\left[Z_1(\alpha x)\right]^2 \quad 18$$

Multiplying 2 and 3 for $x\to\alpha x$, with a little manipulation, recurrence relation is obtained for higher order integrals similar to 17, 18

$$\int x^3\left[Z_{n+1}(\alpha x)\right]^2 dx=\int x^3\left[Z_{n-1}(\alpha x)\right]^2 dx$$
$$+\frac{4n}{\alpha^2}\left\{\frac{x^2}{2}\left[\left[Z_n(\alpha x)\right]^2-Z_{n+1}(\alpha x)Z_{n-1}(\alpha x)\right]\right\}-\frac{2n}{\alpha^2}x^2\left[Z_n(\alpha x)\right]^2 \quad 19$$

From 6,

$$\int\left\{\left[\left(\alpha^2-\beta^2\right)x-\frac{1}{x}\right]Z_1(\alpha x)Z_0(\beta x)\right\}dx$$
$$=\beta xZ_1(\alpha x)Z_{-1}(\beta x)-\alpha xZ_0(\alpha x)Z_0(\beta x)+Z_1(\alpha x)Z_0(\beta x)\equiv Y_{10}^{\alpha\beta}(x) \quad 20$$

Therefore

$$\int xdY_{10}^{\alpha\beta}(x)=\left(\alpha^2-\beta^2\right)\int x^2dxZ_1(\alpha x)Z_0(\beta x)-\int dxZ_1(\alpha x)Z_0(\beta x)$$
$$=xY_{10}^{\alpha\beta}(x)+\int\beta xZ_1(\alpha x)Z_1(\beta x)dx+\int\alpha xZ_0(\alpha x)Z_0(\beta x)dx-\int Z_1(\alpha x)Z_0(\beta x)dx \quad 21$$

This yields the formula

$$W_{10}^{\alpha\beta}(x) \equiv \int x^2 Z_1(\alpha x) Z_0(\beta x)\,dx \;\{\alpha \neq \beta\}$$
$$= -\frac{\beta}{(\alpha^2-\beta^2)} x^2 Z_1(\alpha x) Z_1(\beta x) - \frac{\alpha}{(\alpha^2-\beta^2)} x^2 Z_0(\alpha x) Z_0(\beta x)$$
$$+ \frac{2\alpha^2}{(\alpha^2-\beta^2)^2} x Z_1(\alpha x) Z_0(\beta x) - \frac{2\alpha\beta}{(\alpha^2-\beta^2)^2} x Z_0(\alpha x) Z_1(\beta x) \qquad 22$$

The following result can be derived either from 22 or from 9

$$\int x^2 Z_1(\alpha x) Z_0(\alpha x)\,dx = \frac{1}{2\alpha} x^2 \left[Z_1(\alpha x)\right]^2 \qquad 23$$

The identity 22 can be used to derive the following identities

$$\int x^3 dx \left[Z_0(\alpha x) Z_0(\beta x)\right]$$
$$= \frac{\alpha x^3 Z_1(\alpha x) Z_0(\beta x) - \beta x^3 Z_0(\alpha x) Z_1(\beta x)}{(\alpha^2-\beta^2)}$$
$$+ \frac{4\alpha\beta}{(\alpha^2-\beta^2)^2} x^2 Z_1(\alpha x) Z_1(\beta x) + \frac{2(\alpha^2+\beta^2)}{(\alpha^2-\beta^2)^2} x^2 Z_0(\alpha x) Z_0(\beta x) \qquad 24$$
$$+ \frac{4\beta(\alpha^2+\beta^2)}{(\alpha^2-\beta^2)^3} x Z_0(\alpha x) Z_1(\beta x) - \frac{4\alpha(\alpha^2+\beta^2)}{(\alpha^2-\beta^2)^3} x Z_1(\alpha x) Z_0(\beta x)$$

$$\int x^3 dx \left[Z_1(\alpha x) Z_1(\beta x)\right] = \frac{\beta x^3 Z_0(\alpha x) Z_1(\beta x) - \alpha x^3 Z_1(\alpha x) Z_0(\beta x)}{(\alpha^2-\beta^2)}$$
$$- \frac{4\alpha\beta}{(\alpha^2-\beta^2)^2} x^2 Z_1(\alpha x) Z_1(\beta x) - \frac{2(\alpha^2+\beta^2)}{(\alpha^2-\beta^2)^2} x^2 Z_0(\alpha x) Z_0(\beta x) \qquad 25$$
$$- \frac{4\beta(\alpha^2+\beta^2)}{(\alpha^2-\beta^2)^3} x Z_1(\beta x) Z_0(\alpha x) + \frac{4\alpha(\alpha^2+\beta^2)}{(\alpha^2-\beta^2)^3} x Z_1(\alpha x) Z_0(\beta x)$$

Some of these results have been used in ref. [3].

## III. Identities involving products of three functions:

In the following, functions representing indefinite integrals will be written with explicit functional dependence (x). Define

$$I_{000}^{\alpha\beta\gamma}(x) \equiv \int x Z_0(\alpha x) Z_0(\beta x) Z_0(\gamma x)\,dx \qquad 26$$

$$I_{111}^{\alpha\beta\gamma}(x) \equiv \int xZ_1(\alpha x)Z_1(\beta x)Z_1(\gamma x)dx \quad 27$$

$$I_{110}^{\alpha\beta\gamma}(x) \equiv \int xZ_1(\alpha x)Z_1(\beta x)Z_0(\gamma x)dx \quad 28$$

$$I_{001}^{\alpha\beta\gamma}(x) \equiv \int xZ_0(\alpha x)Z_0(\beta x)Z_1(\gamma x)dx \quad 29$$

$$K_{110}^{\alpha\beta\gamma}(x) \equiv \int Z_1(\alpha x)Z_1(\beta x)Z_0(\gamma x)dx \quad 30$$

$$K_{111}^{\alpha\beta\gamma}(x) \equiv \int Z_1(\alpha x)Z_1(\beta x)Z_1(\gamma x)dx \quad 31$$

$$W_{11}^{\alpha\beta}(x) \equiv \int xdx\left[Z_1(\alpha x)Z_1(\beta x)\right] = \frac{1}{\left(\alpha^2-\beta^2\right)}\left\{\beta xZ_1(\alpha x)Z_0(\beta x) - \alpha xZ_0(\alpha x)Z_1(\beta x)\right\} \quad 32$$

$$W_{00}^{\alpha\beta}(x) \equiv \int xdxZ_0(\alpha x)Z_0(\beta x) = \frac{1}{\left(\alpha^2-\beta^2\right)}\left\{\alpha xZ_1(\alpha x)Z_0(\beta x) - \beta xZ_0(\alpha x)Z_1(\beta x)\right\} \quad 33$$

From 26

$$\begin{aligned} I_{000}^{\alpha\beta\gamma}(x) &\equiv \int Z_0(\alpha x)Z_0(\beta x)xZ_0(\gamma x)dx = \frac{1}{\gamma}\int Z_0(\alpha x)Z_0(\beta x)d\left[xZ_1(\gamma x)\right] \\ &= \frac{1}{\gamma}\left[xZ_1(\gamma x)Z_0(\alpha x)Z_0(\beta x)\right] - \frac{1}{\gamma}\int xZ_1(\gamma x)d\left[Z_0(\alpha x)Z_0(\beta x)\right] \\ &= \frac{1}{\gamma}\left[xZ_1(\gamma x)Z_0(\alpha x)Z_0(\beta x)\right] + \frac{\alpha}{\gamma}I_{110}^{\gamma\alpha\beta}(x) + \frac{\beta}{\gamma}I_{110}^{\beta\gamma\alpha}(x) \end{aligned} \quad 34$$

Therefore,

$$\alpha I_{110}^{\gamma\alpha\beta}(x) + \beta I_{110}^{\beta\gamma\alpha}(x) = \gamma I_{000}^{\alpha\beta\gamma}(x) - \left[xZ_1(\gamma x)Z_0(\alpha x)Z_0(\beta x)\right] \quad 35$$

By cyclic substitution $\alpha \to \beta, \beta \to \gamma, \gamma \to \alpha$ two additional equations are obtained, which can be written in the matrix form

$$\begin{bmatrix} \beta & 0 & \gamma \\ \alpha & \gamma & 0 \\ 0 & \beta & \alpha \end{bmatrix}\begin{bmatrix} I_{110}^{\alpha\beta\gamma}(x) \\ I_{110}^{\beta\gamma\alpha}(x) \\ I_{110}^{\gamma\alpha\beta}(x) \end{bmatrix} = \begin{bmatrix} \alpha I_{000}^{\alpha\beta\gamma}(x) - \left[xZ_1(\alpha x)Z_0(\beta x)Z_0(\gamma x)\right] \\ \beta I_{000}^{\alpha\beta\gamma}(x) - \left[xZ_1(\beta x)Z_0(\gamma x)Z_0(\alpha x)\right] \\ \gamma I_{000}^{\alpha\beta\gamma}(x) - \left[xZ_1(\gamma x)Z_0(\alpha x)Z_0(\beta x)\right] \end{bmatrix} \quad 36$$

This has the solution

$$I_{110}^{\alpha\beta\gamma}(x) = I_{000}^{\alpha\beta\gamma}(x)(2\alpha\beta)^{-1}(\alpha^2+\beta^2-\gamma^2)$$
$$+(2\alpha\beta)^{-1}\left[\gamma x Z_1(\gamma x)Z_0(\alpha x)Z_0(\beta x)-\alpha x Z_1(\alpha x)Z_0(\beta x)Z_0(\gamma x)-\beta x Z_1(\beta x)Z_0(\gamma x)Z_0(\alpha x)\right] \quad 37$$

Two similar expressions can be written by cyclic substitution.

For the special condition $(\alpha^2+\beta^2-\gamma^2)=0$, 37 becomes

$$I_{110}^{\alpha\beta\gamma}(x)\left\{\gamma \to \sqrt{\alpha^2+\beta^2}\right\}$$
$$=(2\alpha\beta)^{-1}\left[\gamma x Z_1(\gamma x)Z_0(\alpha x)Z_0(\beta x)-\alpha x Z_1(\alpha x)Z_0(\beta x)Z_0(\gamma x)-\beta x Z_1(\beta x)Z_0(\gamma x)Z_0(\alpha x)\right] \quad 38$$

An alternate (but not independent) relation between the two integrals can be obtained as follows

$$I_{000}^{\alpha\beta\gamma}(x) = \int x Z_0(\alpha x)Z_0(\beta x)Z_0(\gamma x)dx = \int Z_0(\gamma x)d\left[W_{00}^{\alpha\beta}(x)\right]$$
$$= W_{00}^{\alpha\beta}(x)Z_0(\gamma x) - \int W_{00}^{\alpha\beta}(x)d\left[Z_0(\gamma x)\right] \quad 39$$
$$= W_{00}^{\alpha\beta}(x)Z_0(\gamma x) + \frac{\gamma\alpha}{(\alpha^2-\beta^2)}I_{110}^{\gamma\alpha\beta}(x) - \frac{\beta\gamma}{(\alpha^2-\beta^2)}I_{110}^{\beta\gamma\alpha}(x)$$

$I_{000}^{\alpha\beta\gamma}(x)$ satisfies the following relationship

$$\frac{1}{\beta}\frac{\partial}{\partial\beta}\beta\frac{\partial}{\partial\beta}I_{000}^{\alpha\beta\gamma}(x) = \frac{1}{\alpha}\frac{\partial}{\partial\alpha}\alpha\frac{\partial}{\partial\alpha}I_{000}^{\alpha\beta\gamma}(x) = \frac{1}{\gamma}\frac{\partial}{\partial\gamma}\gamma\frac{\partial}{\partial\gamma}I_{000}^{\alpha\beta\gamma}(x) = -\int x^3 Z_0(\alpha x)Z_0(\beta x)Z_0(\gamma x)dx \quad 40$$

A differential relationship can be demonstrated between the two integrals $I_{000}^{\alpha\beta\gamma}(x)$ and $I_{110}^{\alpha\beta\gamma}(x)$ :

$$\frac{\partial}{\partial\beta}I_{000}^{\alpha\beta\gamma}(x) = -\frac{1}{\alpha}\frac{\partial}{\partial\alpha}\alpha I_{110}^{\alpha\beta\gamma}(x);\ \frac{\partial}{\partial\alpha}I_{000}^{\alpha\beta\gamma}(x) = -\frac{1}{\beta}\frac{\partial}{\partial\beta}\beta I_{110}^{\alpha\beta\gamma}(x) \quad 41$$

This can be used to eliminate one of the integrals giving a differential equation for $I_{000}^{\alpha\beta\gamma}(x)$

$$\frac{\partial}{\partial\alpha} I_{000}^{\alpha\beta\gamma}(x) - F(\alpha,\beta,\gamma) I_{000}^{\alpha\beta\gamma}(x) = G(\alpha,\beta,\gamma,x)$$
$$F(\alpha,\beta,\gamma) \equiv \frac{2\alpha(\beta^2+\gamma^2-\alpha^2)}{(\alpha^2+\beta^2-\gamma^2)^2-4\alpha^2\beta^2} = \frac{\partial}{\partial\alpha}\left\{-\frac{1}{2}\mathrm{Log}\left[(\alpha^2+\beta^2-\gamma^2)^2-4\alpha^2\beta^2\right]\right\}$$
$$G(\alpha,\beta,\gamma,x) \equiv \frac{1}{(\alpha^2+\beta^2-\gamma^2)^2-4\alpha^2\beta^2}$$
$$\times\left\{\begin{array}{l} \gamma(\alpha^2+\beta^2-\gamma^2)x^2 Z_1(\alpha x)Z_0(\beta x)Z_1(\gamma x) \\ +(\alpha^2-\beta^2-\gamma^2)x^2 Z_0(\alpha x)Z_0(\beta x)Z_0(\gamma x) \\ -(\alpha^2+\beta^2-\gamma^2)x^2 Z_1(\alpha x)Z_1(\beta x)Z_0(\gamma x) \\ -2\alpha\beta\gamma x^2 Z_0(\alpha x)Z_1(\beta x)Z_1(\gamma x) \\ +2\alpha^2\beta x^2 Z_1(\alpha x)Z_1(\beta x)Z_0(\gamma x) \end{array}\right\} \quad 42$$

Although the solution of 42 can be formally written down, no simple formula is found for $I_{000}^{\alpha\beta\gamma}(x)$.

Starting from 28,

$$I_{110}^{\alpha\beta\gamma}(x) = \frac{1}{\gamma}\int Z_1(\alpha x)Z_1(\beta x)\,d\left[xZ_1(\gamma x)\right]$$
$$= \frac{1}{\gamma}xZ_1(\gamma x)Z_1(\alpha x)Z_1(\beta x) - \frac{\alpha}{\gamma}\int xZ_1(\gamma x)Z_1(\beta x)\,dxZ_0(\alpha x) \quad 43$$
$$-\frac{\beta}{\gamma}\int xZ_1(\gamma x)Z_1(\alpha x)Z_0(\beta x)\,dx + \frac{2}{\gamma}\int Z_1(\gamma x)Z_1(\beta x)Z_1(\alpha x)\,dx$$

Therefore

$$\gamma I_{110}^{\alpha\beta\gamma}(x) + \alpha I_{110}^{\beta\gamma\alpha}(x) + \beta I_{110}^{\gamma\alpha\beta}(x) = xZ_1(\gamma x)Z_1(\alpha x)Z_1(\beta x) + 2K_{111}^{\alpha\beta\gamma}(x) \quad 44$$

Using 37

$$K_{111}^{\alpha\beta\gamma}(x) = \frac{(\alpha^2-\beta^2-\gamma^2)}{4\beta\gamma}xZ_1(\alpha x)Z_0(\beta x)Z_0(\gamma x) + \frac{(\beta^2-\gamma^2-\alpha^2)}{4\gamma\alpha}xZ_1(\beta x)Z_0(\gamma x)Z_0(\alpha x)$$
$$+\frac{(\gamma^2-\alpha^2-\beta^2)}{4\alpha\beta}xZ_1(\gamma x)Z_0(\alpha x)Z_0(\beta x) - \frac{1}{2}xZ_1(\gamma x)Z_1(\alpha x)Z_1(\beta x) \quad 45$$
$$-I_{000}^{\alpha\beta\gamma}(x)\left\{\frac{(\alpha^2+\beta^2-\gamma^2)^2-4\alpha^2\beta^2}{4\alpha\beta\gamma}\right\}$$

For the special condition $(\alpha^2+\beta^2-\gamma^2)^2-4\alpha^2\beta^2=0$ ,45 becomes

$$K_{111}^{\alpha\beta\gamma}(x)\{\gamma \to |\alpha \pm \beta|\}$$
$$= \frac{(\alpha^2-\beta^2-\gamma^2)}{4\beta\gamma} xZ_1(\alpha x)Z_0(\beta x)Z_0(\gamma x) + \frac{(\beta^2-\gamma^2-\alpha^2)}{4\gamma\alpha} xZ_1(\beta x)Z_0(\gamma x)Z_0(\alpha x)$$
$$+ \frac{(\gamma^2-\alpha^2-\beta^2)}{4\alpha\beta} xZ_1(\gamma x)Z_0(\alpha x)Z_0(\beta x) - \frac{1}{2} xZ_1(\gamma x)Z_1(\alpha x)Z_1(\beta x) \quad 46$$

In a similar manner, the following identities can be proved

$$\gamma I_{001}^{\alpha\beta\gamma}(x) + \alpha I_{001}^{\beta\gamma\alpha}(x) + \beta I_{001}^{\gamma\alpha\beta}(x) = -xZ_0(\alpha x)Z_0(\beta x)Z_0(\gamma x) + K_{000}^{\alpha\beta\gamma} \quad 47$$

$$I_{111}^{\alpha\beta\gamma}(x) = W_{11}^{\alpha\beta}(x)Z_1(\gamma x) - \frac{\beta\gamma}{(\alpha^2-\beta^2)} I_{001}^{\beta\gamma\alpha} + \frac{\gamma\alpha}{(\alpha^2-\beta^2)} I_{001}^{\gamma\alpha\beta} + \frac{\beta}{(\alpha^2-\beta^2)} K_{110}^{\gamma\alpha\beta} - \frac{\alpha}{(\alpha^2-\beta^2)} K_{110}^{\beta\gamma\alpha} \quad 48$$

$$\beta I_{001}^{\beta\gamma\alpha} + \alpha I_{001}^{\gamma\alpha\beta} - K_{110}^{\alpha\beta\gamma} - \gamma I_{111}^{\alpha\beta\gamma}(x) = Z_0(\gamma x)xZ_1(\alpha x)Z_1(\beta x) \quad 49$$

Cyclic substitution in 49 leads to three simultaneous equations which can be used to express $I_{001}^{\alpha\beta\gamma}$ etc in terms of $K_{110}^{\alpha\beta\gamma}$ and $I_{111}^{\alpha\beta\gamma}$ and their cyclic variations.

$$2\alpha\beta I_{001}^{\alpha\beta\gamma}(x) = \alpha K_{110}^{\beta\gamma\alpha}(x) + \beta K_{110}^{\gamma\alpha\beta}(x) - \gamma K_{110}^{\alpha\beta\gamma}(x) + I_{111}^{\alpha\beta\gamma}(x)(\alpha^2+\beta^2-\gamma^2)$$
$$+\alpha xZ_0(\alpha x)Z_1(\beta x)Z_1(\gamma x) + \beta xZ_0(\beta x)Z_1(\gamma x)Z_1(\alpha x) - \gamma xZ_0(\gamma x)Z_1(\alpha x)Z_1(\beta x) \quad 50$$

## IV. Case of Bessel functions of first kind

When the Z functions in the previous section are chosen as Bessel functions of the first kind, definite integrals corresponding to equations 26-33 within the limits {0,1} can be defined. Additionally, α, β, γ are assumed to be zeros of Bessel functions: $\alpha \to j_{q,p}$ which is the $p^{th}$ zero of Bessel function of first kind of order q (q=0 or 1). Define

$${}_qC_{ijk}^{mnp} \equiv \int_0^1 xJ_i(j_{q,m}x)J_j(j_{q,n}x)J_k(j_{q,p}x)dx \quad 51$$

$${}_qD_{ijk}^{mnp} \equiv \int_0^1 J_i(j_{q,m}x)J_j(j_{q,n}x)J_k(j_{q,p}x)dx \quad 52$$

From 37

$$ {}_{q}C_{110}^{mnp} = {}_{q}C_{000}^{mnp}\,\frac{\left(j_{q,m}^2 + j_{q,n}^2 - j_{q,p}^2\right)}{2 j_{q,m} j_{q,n}} \qquad 53 $$

From 45

$$ {}_{q}D_{111}^{mnp} = -{}_{q}C_{000}^{mnp}\left\{\frac{\left(j_{q,m}{}^2 + j_{q,n}{}^2 - j_{q,p}{}^2\right)^2 - 4 j_{q,m}{}^2 j_{q,n}{}^2}{4 j_{q,m} j_{q,n} j_{q,p}}\right\} \qquad 54 $$

No formula has yet been found for ${}_{q}C_{000}^{mnp}$. An approximate scaling relations for ${}_{1}C_{000}^{mnp}$ is obtained in the next section.

The importance of these integrals may be appreciated from the following. The product of two Bessel functions may be expanded into a Fourier-Bessel series:

$$ J_j\left(j_{i,n}x\right) J_k\left(j_{i,p}x\right) = \sum_{p=1}^{\infty} c_{ijk}^{mnp} J_i\left(j_{i,p}x\right) \qquad 55 $$

The coefficients of this expansion are given by

$$ c_{ijk}^{mnp} = 2\left[J_{i+1}\left(j_{i,p}\right)\right]^{-2}{}_{i}C_{ijk}^{mnp} \qquad 56 $$

As an illustration, consider the generalized Ohm’s law for a two-fluid plasma with electron inertia [5]

$$ \vec{E} + \vec{v}\times B = \eta J + \frac{}{en}\left\{ J\times B - \nabla p_e - \nabla\cdot\left\{\frac{}{\varepsilon_0\omega_p^2}\right\}\right\} + \frac{}{\varepsilon_0\omega_p^2}\left\{\frac{}{\partial t} + \nabla\cdot\left(J\vec{v} + \vec{v}J\right)\right\} \qquad 57 $$

Using Maxwell’s equations, this can be written as

$$ \frac{\partial^2\vec{E}}{\partial t^2} + \nu\frac{\partial\vec{E}}{\partial t} + \omega_p^2\vec{E} + \omega_p^2\vec{v}\times B + \nu c\,\vec{\nabla}\times\vec{B} + c^2\vec{\nabla}\times\vec{\nabla}\times\vec{E} - \frac{}{\varepsilon_0 m_e}\vec{\nabla}p_e $$
$$ = \frac{e}{\varepsilon_0 m_e}\vec{J}\times\vec{B} + \varepsilon_o^{-1}\vec{\nabla}\cdot\left(\vec{J}\vec{v} + \vec{v}J\right) - \left\{\frac{}{en\varepsilon_0}\right\}\vec{\nabla}\cdot\left(\vec{J}\vec{J}\right) + \left\{\frac{}{en^2\varepsilon_0}\right\}\vec{J}\left(\vec{J}\cdot\vec{\nabla}n\right) \qquad 58 $$

Expressing the electric and magnetic fields in terms of C-K functions in cylindrical geometry with azimuthal symmetry having time dependent mode amplitudes, a system of time evolution equations can be obtained whose left hand side consists of a linear differential operator acting on the mode amplitudes and right hand side consists of a “sum over modes” of terms involving products of quantities related to the mode amplitudes weighted by coefficients of the form given by 51 and 52. In the case being considered, it turns out that all the coefficients are of one of the three form: ${}_{1}C_{000}^{mnp}$ or ${}_{1}C_{110}^{mnp}$ or ${}_{1}D_{111}^{mnp}$.

This procedure allows taking into account effects due to the electron momentum convection term $\vec{\nabla}\cdot\left(\vec{J}\vec{J}\right)$ which is of order $c^2\omega_p^{-2}\ell^{-}$ as compared with the Hall effect term, $\ell$ being a gradient scale length and $\omega_p$ is the electron plasma frequency. This term is usually too stiff to be used in the conventional approach using partial differential equations. The difficulty in numerical treatment of this term in a system of difference equation is converted into the difficulty of evaluating the coefficients of the form ${}_1C_{000}^{mnp}$ or ${}_1C_{110}^{mnp}$ or ${}_1D_{111}^{mnp}$ for large mode numbers (finer spatial scale). Identities 53 and 54 reduce the problem to evaluating only ${}_1C_{000}^{mnp}$ as a function of mode numbers. The electron momentum convection term, being a gradient, should give a higher contribution at higher mode numbers. Since the sum over modes has to be evaluated for a very large number of modes (~$10^7$), a relationship showing how ${}_1C_{000}^{mnp}$ scales with the mode numbers m, n, p is very desirable to enable estimation of the net effect of the electron momentum convection term and similar other terms to observable average effects.

## V. Approximate scaling relation for ${}_1C_{000}^{pp'p''}$

In the absence of an analytic formula, an approximate scaling relation is derived below for ${}_1C_{000}^{mnp}$ for the case of large mode numbers, where the integrand oscillates very rapidly creating difficulty with numerical evaluation of integrals. For smaller mode numbers, the coefficients can be numerically evaluated and tabulated once for all.

The starting point is the result that p<sup>th</sup> zero of Bessel function of first kind of order 1 is $j_p$~$p\pi$ for large p and the asymptotic approximation for Bessel functions of first kind for integral orders [6] for large values of the argument:

$$J_n(x) \approx \sqrt{\frac{2}{\pi x}}\left[\cos\left(x-\tfrac{1}{2}n\pi-\tfrac{1}{4}\pi\right)\right] \qquad \text{for large x} \qquad 59$$

This leads to the following useful approximate formula

$$\Im_{pp'} \equiv \int_0^1 \xi^2 J_0\left(j_{p'}\xi\right) J_0\left(j_p x\right) d\xi = \approx \begin{cases} \dfrac{2}{\pi\sqrt{j_{p'} j_p}}\left\{ -\dfrac{1}{\left(j_{p'} - j_p\right)^2} + \dfrac{(-1)^{p+p'}}{\left(j_{p'} - j_p\right)^2} - \dfrac{(-1)^{p+p'}}{\left(j_{p'} + j_p\right)} \right\} & \text{for } p \neq p' \\ \dfrac{1}{2\pi j_p} & \text{for } p = p' \end{cases} \qquad 60$$

The triple product integral can be approximated as follows:

$$C_{ijk}^{mnp} \equiv \int_0^1 \xi d\xi J_i\left(j_m\xi\right) J_j\left(j_n\xi\right) J_k\left(j_p\xi\right)$$

$$\approx \frac{1}{4}\frac{2\sqrt{2}}{\pi^3\sqrt{mnp}} \int_0^1 \frac{d\xi}{\sqrt{\xi}} \begin{Bmatrix} \cos\left[\pi\left(m+n+p\right)\xi - \tfrac{1}{2}\left(i+j+k\right)\pi - \tfrac{3}{4}\pi\right] \\ +\cos\left[\pi\left(m-n-p\right)\xi - \tfrac{1}{2}\left(i-j-k\right)\pi + \tfrac{1}{4}\pi\right] \\ +\cos\left[\pi\left(m+n-p\right)\xi - \tfrac{1}{2}\left(i+j-k\right)\pi - \tfrac{1}{4}\pi\right] \\ +\cos\left[\pi\left(m-n+p\right)\xi - \tfrac{1}{2}\left(i-j+k\right)\pi - \tfrac{1}{4}\pi\right] \end{Bmatrix} \qquad 61$$

It can be shown that

$$\int_0^1 \frac{d\xi}{\sqrt{\xi}} \cos\left(P\pi\xi + Q\pi\right)$$

$$= 2\cos\left(Q\pi\right)\left\{ \delta_{P,0} + \left(1-\delta_{P,0}\right)\frac{C\left(\sqrt{2|P|}\right)}{\sqrt{2|P|}} \right\} - 2\sin\left(Q\pi\right)\text{sign}\left(P\right)\frac{S\left(\sqrt{2|P|}\right)}{\sqrt{2|P|}} \qquad 62$$

where the Fresnel Integrals are defined by

$$S\left(t\right) = \int_0^t \sin\left(\pi t'^2/2\right) dt'; \quad C\left(t\right) = \int_0^t \cos\left(\pi t'^2/2\right) dt'$$

Applying this to 61

$$
C_{ijk}^{mnp} = \frac{1}{2\pi^3\sqrt{mnp}} \left\{ \begin{array}{l}
\cos\left(-\tfrac{1}{2}(i+j+k)\pi - \tfrac{3}{4}\pi\right)\left\{\frac{C\left(\sqrt{2\left|(m+n+p)\right|}\right)}{\sqrt{\left|(m+n+p)\right|}}\right\} \\
-\frac{\sin\left(-\tfrac{1}{2}(i+j+k)\pi - \tfrac{3}{4}\pi\right)S\left(\sqrt{2\left|(m+n+p)\right|}\right)}{\sqrt{\left|(m+n+p)\right|}} \\
+\cos\left(-\tfrac{1}{2}(i-j-k)\pi + \tfrac{1}{4}\pi\right)\left\{\delta_{(m-n-p),0} + \left(1-\delta_{(m-n-p),0}\right)\frac{C\left(\sqrt{2\left|(m-n-p)\right|}\right)}{\sqrt{\left|(m-n-p)\right|}}\right\} \\
-\mathrm{sign}(m-n-p)\frac{\sin\left(-\tfrac{1}{2}(i-j-k)\pi + \tfrac{1}{4}\pi\right)S\left(\sqrt{2\left|(m-n-p)\right|}\right)}{\sqrt{\left|(m-n-p)\right|}} \\
+\cos\left(-\tfrac{1}{2}(i+j-k)\pi - \tfrac{1}{4}\pi\right)\left\{\delta_{(m+n-p),0} + \left(1-\delta_{(m+n-p),0}\right)\frac{C\left(\sqrt{2\left|(m+n-p)\right|}\right)}{\sqrt{\left|(m+n-p)\right|}}\right\} \\
-\mathrm{sign}(m+n-p)\frac{\sin\left(-\tfrac{1}{2}(i+j-k)\pi - \tfrac{1}{4}\pi\right)S\left(\sqrt{2\left|(m+n-p)\right|}\right)}{\sqrt{\left|(m+n-p)\right|}} \\
+\cos\left(-\tfrac{1}{2}(i-j+k)\pi - \tfrac{1}{4}\pi\right)\left\{\delta_{(m-n+p),0} + \left(1-\delta_{(m-n+p),0}\right)\frac{C\left(\sqrt{2\left|(m-n+p)\right|}\right)}{\sqrt{\left|(m-n+p)\right|}}\right\} \\
-\mathrm{sign}(m-n+p)\frac{\sin\left(-\tfrac{1}{2}(i-j+k)\pi - \tfrac{1}{4}\pi\right)S\left(\sqrt{2\left|(m-n+p)\right|}\right)}{\sqrt{\left|(m-n+p)\right|}}
\end{array} \right\} \quad 63
$$

## VI. Numerical Validation

The identities derived above have been validated by computing the left and right hand sides separately for varying arguments. The difference in the LHS and RHS ranges from $10^{-8}$ to $10^{-22}$; the poorer correspondence is correlated with the divergence of the Bessel function of second kind at small values of the argument. The approximate nature of relation 60 and numerical evaluation of the derivative in 41 limit the agreement to within $10^{-4}$-$10^{-5}$. The most difficult challenge is posed by validation of 63, for the case of ${}_1C_{000}^{mnp}$. The approximation represented by 59 is expected to be more accurate for larger mode numbers while the numerical evaluation of $C_{ijk}^{mnp}$ becomes very difficult for large mode numbers. As a compromise, numerically integrated values and those calculated from the approximation formula are compared in Table 1 and fig. 1 for various combinations of the mode indices. Above mode number 150, the numerical integration converges very slowly and for some mode combinations, fails the convergence test. The agreement between the numerical integration and formula 63 is seen to be better at higher mode numbers than at lower mode numbers, as expected. However, even at higher mode numbers, there is some disagreement, when the numerical integration gives a negative number while the formula doesn't. It is not known whether this is due to numerical error in the integration because of a highly oscillatory integrand or because of the error in the approximation. The inclusion of higher order terms in 59 is not expected to improve matters very much, firstly because their analytical integration in 61 is not straightforward and secondly because their contribution at higher mode numbers is expected to diminish faster than that of the dominant term.

**Table 1:**

| TABLE-1 | | | | |
|---|---|---|---|---|
| m | n | p | LHS | RHS |
| 44 | 23 | 63 | 4.557E-05 | 3.140E-05 |
| 20 | 20 | 20 | 9.061E-05 | 8.071E-05 |
| 22 | 22 | 22 | 7.508E-05 | 6.689E-05 |
| 25 | 25 | 30 | 5.265E-05 | 4.659E-05 |
| 30 | 30 | 30 | 4.065E-05 | 3.627E-05 |
| 22 | 89 | 31 | -7.053E-10 | 9.828E-06 |
| 40 | 45 | 50 | 1.866E-05 | 1.787E-05 |
| 37 | 77 | 57 | 1.590E-05 | 1.408E-05 |
| 47 | 61 | 87 | 1.162E-05 | 1.032E-05 |
| 29 | 47 | 57 | 2.342E-05 | 2.158E-05 |
| 40 | 40 | 40 | 2.297E-05 | 2.053E-05 |
| 50 | 50 | 50 | 1.474E-05 | 1.320E-05 |
| 60 | 60 | 60 | 1.025E-05 | 9.195E-06 |
| 70 | 70 | 70 | 7.543E-06 | 6.770E-06 |
| 80 | 80 | 80 | 5.781E-06 | 5.195E-06 |
| 90 | 90 | 90 | 4.571E-06 | 4.112E-06 |
| 100 | 100 | 100 | 3.704E-06 | 3.335E-06 |
| 105 | 85 | 97 | 4.139E-06 | 3.899E-06 |
| 100 | 50 | 20 | -1.114E-09 | 7.515E-06 |
| 120 | 120 | 120 | 2.575E-06 | 2.321E-06 |
| 130 | 130 | 130 | 2.195E-06 | 1.980E-06 |
| 140 | 140 | 140 | 1.893E-06 | 1.708E-06 |
| 150 | 150 | 150 | 1.649E-06 | 1.489E-06 |
| 160 | 160 | 160 | 1.450E-06 | 1.310E-06 |
| 160 | 80 | 70 | -1.052E-08 | 2.131E-06 |
| 170 | 170 | 170 | 1.285E-06 | 1.161E-06 |
| 200 | 200 | 185 | 9.807E-07 | 9.184E-07 |
| 200 | 200 | 200 | 9.286E-07 | 8.401E-07 |
| 200 | 100 | 185 | 1.749E-06 | 1.582E-06 |

Comparison of LHS and RHS of 63. The LHS is calculated by numerical integration over the product of three Bessel functions. RHS is the approximation formula. The trial values of mode numbers m,n,p are chosen so as to get a feel of how the integral behaves in various situations.

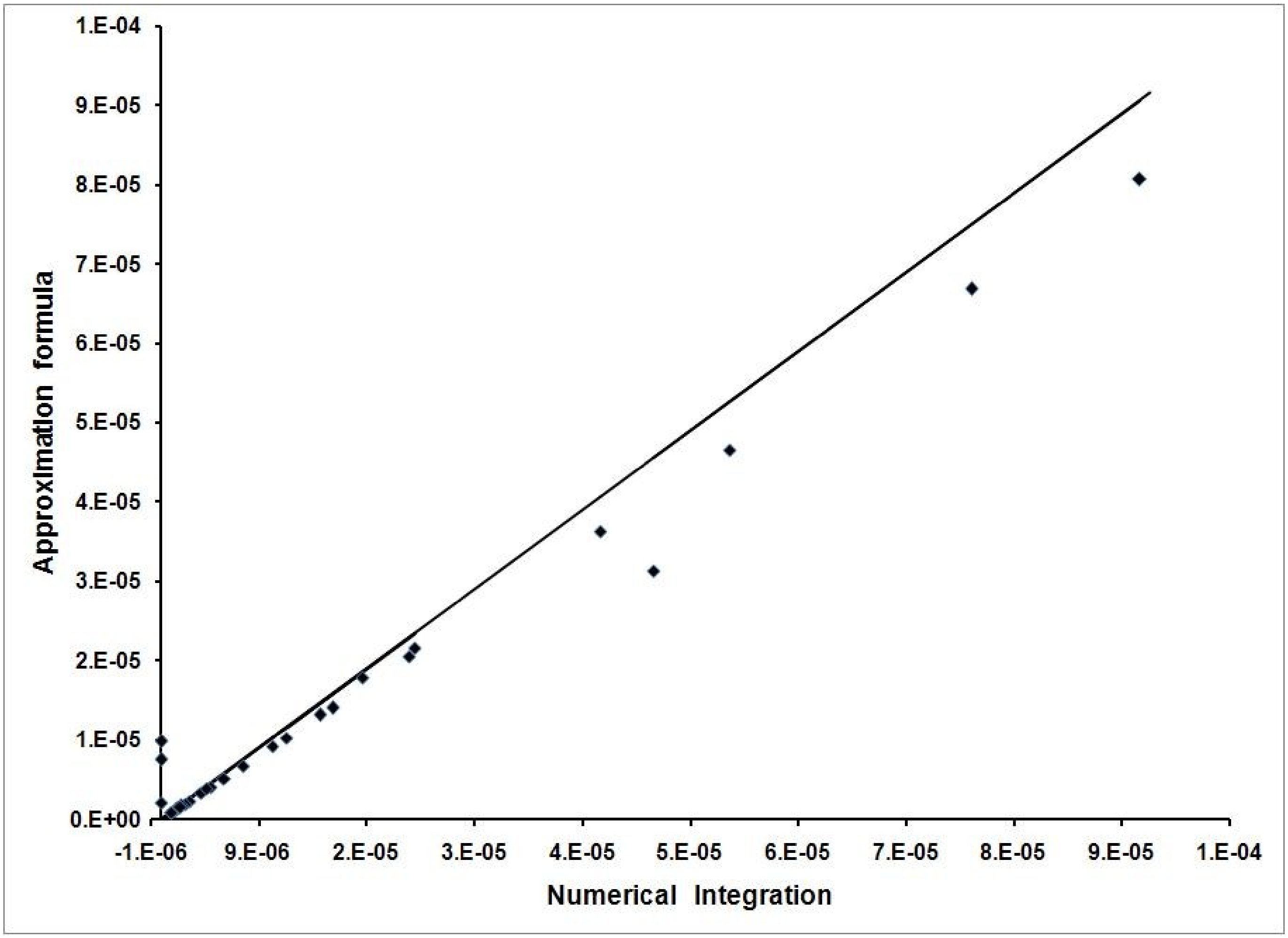

Fig-1

Comparison of the approximation formula with the numerical integration. The trend line has slope 1. The approximation formula tends to underestimate the integral. The three points nearly on the y-axis actually have a value close to zero and the numerical integration shows errors in convergence.